\documentclass[11pt]{article}
\pdfoutput=1
\usepackage{caption}
\usepackage{url}
\usepackage{subfigure}
\usepackage[T1]{fontenc}
\usepackage[utf8]{inputenc}
\usepackage{authblk}
\usepackage{color}
\usepackage{amsmath, amssymb, latexsym}
\usepackage{psfrag,epsfig,amsfonts,amsmath,latexsym,amsthm,amssymb,amscd,url }
\usepackage[colorlinks=true,linkcolor=blue,citecolor=blue, linktocpage=true]{hyperref}

\renewcommand{\phi}{\varphi}

\newcommand{\BE}{\begin{equation}}
\newcommand{\EE}{\end{equation}}
\newcommand{\BEN}{\begin{equation*}}
\newcommand{\EEN}{\end{equation*}}
\newcommand{\BAL}{\begin{align}}
\newcommand{\EAL}{\end{align}}
\newcommand{\BAN}{\begin{align*}}

\begin{document}

\title{Extracting Stochastic Governing Laws by Nonlocal Kramers-Moyal Formulas}

\author[a]{Yubin Lu}
\author[b]{Yang Li}
\author[c]{Jinqiao Duan \thanks{Corresponding author: duan@iit.edu}}

\affil[a]{School of Mathematics and Statistics \& Center for Mathematical Sciences, Huazhong University of Science and Technology, Wuhan 430074, China}
\affil[b]{School of Automation, Nanjing University of Science and Technology, Nanjing 210094, China}
\affil[c]{Departments of Applied Mathematics \& Physics, Illinois Institute of Technology, Chicago, IL 60616, USA}

\renewcommand*{\Affilfont}{\small\it}
\renewcommand\Authands{ and }
\date{\today}
\maketitle

 \begin{abstract}
With the rapid development of computational techniques and scientific tools, great progress of data-driven analysis has been made to extract governing laws of dynamical systems from data. Despite the wide occurrences of non-Gaussian fluctuations, the effective data-driven methods to identify stochastic differential equations with non-Gaussian L\'evy noise are relatively few so far. In this work, we propose a data-driven approach to extract stochastic governing laws with both (Gaussian) Brownian motion and (non-Gaussian) L\'evy motion, from short bursts of simulation data. Specifically, we use the normalizing flows technology to estimate the transition probability density function (solution of nonlocal Fokker-Planck equation) from data, and then substitute it into the recently proposed nonlocal Kramers-Moyal formulas to approximate L\'evy jump measure, drift coefficient and diffusion coefficient. We demonstrate that this approach can learn the stochastic differential equation with L\'evy motion. We present examples with one- and two-dimensional, decoupled and coupled systems to illustrate our method. This approach will become an effective tool for discovering stochastic governing laws and understanding complex dynamical behaviors.\\
\end{abstract}

\bigskip\par



\section{Introduction}

Differential equations are usually used to model complex phenomena in climate, molecular biology, condensed matter physics, mechanical systems and other fields. In order to improve the understanding and make predictions about the future dynamical evolution, researchers need to describe the effects of noise in dynamical systems. Stochastic differential equations are an effective formalism for modeling physical phenomena under random fluctuations \cite{Duan2015}. It is often assumed that the random fluctuations of deterministic dynamical systems are Gaussian   noise. In the last decade or two, however, researchers began to actively study the dynamical behaviors for systems under   non-Gaussian influences.

Stochastic differential equations with non-Gaussian L\'evy motions have the potential for wider applications in physical sciences. For instance, stochastic dynamical systems with both (Gaussian) Brownian motion and (non-Gaussian) L\'evy motion can be used to model a significant and widely used class of Markov processes, so-called Feller processes \cite{Bottcher}. According to the Greenland ice core measurement data, Ditlevsen \cite{Ditlevsen} found that the temperature in that climate system could be modeled as stochastic differential equations with $\alpha$-stable L\'evy motion. Other researchers also used the non-Gaussian L\'evy motion to characterize the random fluctuations emerged in neural systems \cite{NeuralModel}, gene networks \cite{GeneModel}, and the Earth systems \cite{ashwin2012tipping,YangDuanWiggins2020,ZhengYY2020,WeiPY2020,DannyTesfay}.

However, it is often difficult to establish the mathematical models of complex systems in nature, especially when less-understood or unknown mechanisms are present. Thanks to recent improvements in the storage and computing power of computers, more experimental, observational and simulated data are now available. Data-driven approaches to complex dynamical systems have attracted increasing attention in recent years, including but not limited to system identification \cite{maulik2020time,SINDy1,SINDy2,GPs1,GPs2,DaiMinChaos,YangLi2020a,LuYB2020,DMD,EDMD,SKO1,SKO2,SKO3,Kevrekidis1,Kevrekidis4}, model reduction \cite{FeiLu,SKO3,Klus1,Klus2,Kevrekidis2,Kevrekidis8,Kevrekidis10,ChenBH2009}, extracting dynamical behaviors \cite{Froyland,Froyland2009,Koltai2017,Dellnitz,Metzner,Tantet,Thiede,Kevrekidis6} directly from data.

System identification of differential equations models for physical phenomena is a major task in data-driven analysis of complex dynamical systems. For deterministic systems, Brunton et. al \cite{SINDy1,SINDy2} developed tools to identify  models as ordinary differential equations. For stochastic     systems, there are at least three common ways to learn or extract the coefficients of  stochastic differential equations models. The first way is to directly estimate these coefficients where mainstream methods include (nonlocal)  Kramers-Moyal formulas \cite{DaiMinChaos,YangLi2020a,KM} , neural differential equations \cite{Duvenaud1,Duvenaud2,NeuralSDE,Felix} and Bayesian inference \cite{GPs1,GPs2}. The second way is to learn the generator of stochastic differential equation by estimating the stochastic Koopman operator \cite{LuYB2020,SKO1,SKO2,SKO3}. Subsequently, we can obtain the coefficients of stochastic differential equation. In the third way, the coefficients are obtained by estimating the (nonlocal) Fokker-Planck equations corresponding to stochastic differential equations. For example,     Chen et. al and Yang et. al  \cite{XiaoliChen,LiuYang} identify   stochastic systems by estimating the coefficients of the associated nonlocal Fokker-Planck equations.

  The usual Kramers-Moyal formulas link the sample path data with coefficients of the stochastic systems with Gaussian noise. Therefore, it is possible to identify the systems by using the sample path data of stochastic differential equations. For stochastic dynamical systems with non-Gaussian fluctuations, We devised the so-called nonlocal Kramers-Moyal formulas   in our previous work \cite{YangLi2020a}. Hence the nonlocal Kramers-Moyal formulas can be used to identify stochastic differential equations with non-Gaussian fluctuations (in particular, L\'evy motions). However, the original algorithm of nonlocal Kramers-Moyal formulas required a lot of sample data and assumptions that imposed a basis on the drift and diffusion coefficients.

  In this present work, we devise a method to identify stochastic dynamical systems aiming to reduce data requirements and discard the basis assumption of the coefficients. Concretely, we first estimate the probability density (via normalizing flow neural networks; see below) of a stochastic differential equation with L\'evy motion. Subsequently, we apply the nonlocal Kramers-Moyal formulas \cite{YangLi2020a}    to calculate the L\'evy jump measure, drift coefficient and diffusion coefficient of the stochastic differential equation.

For this new method, probability density estimation plays an important role and this can be seen from nonlocal Kramers-Moyal formulas \cite{YangLi2020a}. In this present work, we are inspired by modern \emph{generative} machine learning techniques that allow for the ability to learn the underlying distribution that generates samples of a data set. Examples of such methods include variational autoencoders (VAEs), generative adversarial networks (GANs) and normalizing flows \cite{NFs1,NFs2,IAF,MAF,RealNVP,stNFs,LuYB2021}, etc. For VAEs and GANs, however, they can not provide an explicit density estimation. But the explicit density estimation is exactly what we need. Fortunately, normalizing flow is exactly what we want.

The idea of a normalizing flow was introduced by Tabak and Vanden-Eijnden \cite{tabak-2010} in the context of density estimation and sampling. A flow that is "normalizing" transforms a given
distribution, generally into the standard normal one, while its (invertible) Jacobian is traced and inverted continuously to obtain the inverse map. By constructing the transformation that meets our requirements, we can get an estimate of the target probability density.

To demonstrate the performance of our method, we illustrate     with one- and two-dimensional stochastic differential equations with L\'evy motions. Furthermore, we also consider the effects of different types of noise, i.e., additive noise and multiplicative noise.

The remainder of this article is arranged as follows: In Section 2, we introduce the idea of normalizing flows and provide two types of transformations to show that  how to estimate probability density from sample data. In Section 3, we review the nonlocal Kramers-Moyal formulas, in order to compute the L\'evy jump measure, drift and diffusion by estimated density. In Section 4, we will show our results for different cases and compare with true coefficients. The final section contains a discussion of this study and perspectives for building on this line of research. \\

\section{Density estimation }
We recall normalizing flows in subsection \ref{NF}, and then introduce two concrete cases of normalizing flows. One is neural spline flows in subsection \ref{NSF}, the other one is real-value non-volume preserving (RealNVP) methodology in the subsection \ref{RealNVP}. For those examples shown in Section \ref{results}, the neural spline flows are used to estimate one-dimensional probability density and the RealNVP are used to estimate two-dimensional probability densities.
\subsection{Normalizing flows}\label{NF}
Normalizing flows is a general idea to express probability density  using a prior probability density and a series of bijective transformations. In the following, we adopt the notations from \cite{NFs2}. Given $z$ as a $D-$dimensional real vector sampled from $p_{z}(z)$. The main idea of normalizing flow is to find a transformation $T$ such that:
\begin{align}
z=T(x), \quad where \quad z\thicksim p_{z}(z),
\end{align}
where $p_{z}(z)$ is a prior (or latent) density. When the transformation $T$ is invertible and both $T$ and $T^{-1}$ are differentiable, the density $p_{x}(x)$ can be calculated by a change of variables:
\begin{align}\label{CoV}
p_{x}(x)=p_{z}(T(x))\mid {\rm det} J_{T}(x)\mid, \quad where \quad z=T(x)
\end{align}
and $J_{T}(x)$ is the Jacobian of $T$, i.e.,

\begin{center} {$J_{T}(x) = \left[ {\begin{array}{*{20}{c}}
 \frac{\partial T_1}{\partial x_1} & \cdots & \frac{\partial T_1}{\partial x_D}\\
 \vdots & \ddots & \vdots\\
 \frac{\partial T_D}{\partial x_1} & \cdots & \frac{\partial T_D}{\partial x_D}
\end{array}} \right].$}
\end{center}

For two invertible and differentiable transformations $T_1$ and $T_2$, we have the following properties:
\begin{align}
(T_{2}\circ T_{1})^{-1} &= T_{1}^{-1}\circ T_{2}^{-1} \nonumber\\
{\rm det} J_{T_{2}\circ T_{1}(z)}&={\rm det} J_{T_{2}}(T_{1}(z))\cdot {\rm det} J_{T_{1}}(z). \nonumber
\end{align}
Consequently, we can construct complex transformations by composing multiple instance of simpler transformations, i.e., $T=T_{K}\circ T_{K-1}\circ \cdots \circ T_{1}$, where each $T_K$ transforms $z_{K-1}$ into $z_K$, assuming $z_0=x$ and $z_K=z$.

Let the set of samples $\{x_0, x_1, \ldots, x_n \}$  be taken from an unknown distribution $p_{x}(x;\theta)$, we can minimize the negative log-likelihood on data,
\begin{align}\label{loss1}
\mathcal{L} = -\sum_{i=1}^n {\rm log} p_{x}(x_i; \theta).
\end{align}
Take the change of variables formula (\ref{CoV}) into (\ref{loss1}), we have
\begin{align}\label{loss2}
\mathcal{L} = -\sum_{i=1}^n [{\rm log} p_{z}(T(x_i)) + {\rm log}\mid {\rm det} J_{T}(x)\mid_{x=x_i}].
\end{align}
Therefore, we can learn the transformation $T$ by minimizing the loss function (\ref{loss2}) .

This subsection \ref{NF} introduced the general idea of normalizing flows. We can see that the transformation $T$ plays a crucial role in approximating the target density function. Now we introduce two specific transformations, i.e., neural spline flows and RealNVP.

\subsection{Neural spline flows}\label{NSF}
In this subsection we briefly review the neural spline flows. For the sake of simplicity, we only consider one-dimensional density estimation here. The authors Durkan et al \cite{NSF} proposed to implement the transformation $T$ using monotonic rational-quadratic splines. A rational-quadratic function takes the form of a quotient of two quadratic polynomials.

We adopt the notations from \cite{NSF}. The spline uses $K$ different rational-quadratic functions, with boundaries set by $K+1$ coordinates $\{x^{(k)}, y^{(k)}\}_{k=0}^K$ known as knots. The knots monotonically increase between $(x^{(0)}, y^{(0)})=(-B, -B)$ and $(x^{(K)}, y^{(K)})=(B, B)$. The spline itself maps an interval $[-B, B]$ to $[-B, B]$. The $K-1$ derivatives at the internal points are set to arbitrary positive values and set the boundaries derivatives to $1$.

In short, we can obtain a unique spline if we give the widths $\theta^w$ , heights $\theta^h$ of the $K$ bins and the $K-1$ derivatives $\theta^d$ at the internal knots.

Denoting $\theta=[\theta^w,\theta^h,\theta^d]$, the author X constructed the transformation as follows,
\begin{align}\label{NSF_trans}
   \theta &= NN(x) \nonumber\\
    z&=T_{\theta}(x),
\end{align}
where $NN$ is a neural network and $T_{\theta}$ is a monotonic rational-quadratic spline.

Therefore, we can estimate the target density $p_{x}(x)$ by the rational-quadratic spline $T_{\theta}$ and the prior density $p_{z}(z)$
\begin{align}
    p_{x}(x)&=p_{z}(T(x))\mid {\rm det} J_{T_{\theta}}(x)\mid.
\end{align}
For more details, see \cite{NSF}.

\subsection{Real-value non-volume preserving transformations}\label{RealNVP}
After briefly reviewing the general framework of normalizing flows and neural spline flows in previous subsections, we introduce the construction of the real-value non-volume preserving transformation \cite{RealNVP}. In this subsection, for the sake of simplicity, we only focus on two-dimensional density estimation. In fact, this method can be extended to any dimension.\\

For the target probability density $p_{x}(x)$ , where $x=(x_1, x_2)$ is a two-dimensional real vector and the prior density is denoted by $p_{z}(z)$, $z\in \mathbb{R}^2$, we aim to design an invertible and differentiable transformation $T_{\theta}$ such that $z=T_{\theta}(x)$ and determinant of the Jacobian is easy to compute. To be more specific, we propose the following transformation,
\begin{align}\label{RealNVP_Trans}
z_{1}&=x_{1}, \nonumber \\
z_{2}&=\frac{1}{C}[x_{2} e^{\mu(x_{1})}+\nu(x_{1})],
\end{align}
where the notation $\mu$ and $\nu$ are two different neural networks. Here $C$ is a hyperparameter. We denote this transformation by $T_\theta$. The Jacobian matrix then becomes a computationally tractable
\begin{align}
    {\rm det}J_{T_{\theta}} = e^{\mu(x_{1},t)}.
\end{align}

In consequence, for a given dataset  $\mathcal{D}=\{x_{i}\}_{i=1}^{n}$ sampled from $p_{x}(x)$, the negative log-likelihood becomes
\begin{align}\label{loss3}
\mathcal{L} &= -\sum_{i=1}^n {\rm log} p_{x}(x_{i}; \theta) \nonumber \\
&=-\sum_{i=1}^n  [{\rm log} p_{z}(T_{\theta}(x_{i})) + {\rm log}\mid {\rm det} J_{T_{\theta}}(x)\mid_{x=x_{i}}].
\end{align}
Minimizing the loss function (\ref{loss3}) to get the optimal parameters of neural networks $\mu$ and $\nu$, i.e., we learn a flexible, yet expressive, transformation $T_{\theta}$ from dataset $\mathcal{D}$. Therefore, we can use this transformation $x=T_{\theta}^{-1}(z)$ to resample from  probability density $p_x(x)$. This means that we have learned a generative model that mimics the underlying random variable that has generated our samples in the training data set.

So far, we have seen how we use normalizing flows to estimate the probability density. Next, we will introduce the nonlocal Kramers-Moyal formulas and explain how to combine probability density with extracting stochastic governing laws.\\

\section{Nonlocal Kramers-Moyal formulas}
\label{KMformulas}
Last section introduces normalizing flows devoted to approximating the probability density function from sample data. It still requires to compute the L\'evy jump measure, drift and diffusion terms based on the estimated probability density, in order to extract the stochastic dynamical systems. And this step can be accomplished by the recently proposed nonlocal Kramers-Moyal formulas. We will review it in the following.

Random fluctuations often have both Gaussian and non-Gaussian statistical features. By L\'evy-It\^o decomposition theorem \cite{Duan2015, Applebaum}, a large class of random fluctuations are indeed modeled as linear combinations of a (Gaussian) Brownian motion $B_t$ and a (non-Gaussian) L\'evy process $L_t$. We thus consider an $n$-dimensional stochastic dynamical system in the following form
\begin{equation} \label{SDE}
dx\left( t \right)=b\left(x\left( t \right) \right)dt+\Lambda \left( x\left( t \right) \right)d{{B}_{t}}+\sigma d{{L}_{t}},
\end{equation}
where ${{B}_{t}}={{\left[ {{B}_{1,t}},\ \cdots ,\ {{B}_{n,t}} \right]}^{T}}$  is an $n$-dimensional Brownian motion,  and ${{L}_{t}}={{\left[ {{L}_{1,t}},\ \cdots ,\ {{L}_{n,t}} \right]}^{T}}$ is an $n$-dimensional non-Gaussian L\'evy process with independent components described in the Appendix. The vector $b\left(x \right)={{\left[ {{b}_{1}}\left( x \right),\ \cdots ,\ {{b}_{n}}\left( x \right) \right]}^{T}}$  is the drift coefficient (or vector field) in ${{\mathbb{R}}^{n}}$  and $\Lambda \left( x \right)$ is an $n\times n$  matrix.  The diffusion matrix is defined by $a\left( x \right)=\Lambda {{\Lambda }^{T}}$.  We take the positive constant $\sigma $  as  the noise intensity of the L\'evy process. Assume that the initial condition is $x\left( 0 \right)=z$  and the jump measure of ${{L}_{t}}$ is ${{\nu }}\left( dy \right)={{W}}\left( y \right)dy$, with kernel $W$,  for $y\in \mathbb{R}^{n}\backslash \left\{ 0 \right\}$.

Based on the Fokker-Planck equation for the probability density function of the system (\ref{SDE}), Li and Duan \cite{YangLi2020a} derived the nonlocal Kramers-Moyal formulas to express the L\'evy jump measure, drift coefficient and diffusion coefficient via sample path. Its kernel idea is to split the phase space into the regions of big and small jumps by a spherical surface. The integrals inside the ball (instead of the whole space in Kramers-Moyal formulas) provide the drift and diffusion, and the jump measure is computed by the data outside, respectively. We restate their main results including the following theorem and corollary as the theoretical foundation of our method.

The following theorem aims to express the L\'evy jump measure, drift and diffusion in terms of the solution of Fokker-Planck equation $p\left( x,t|z,0 \right)$.

\newtheorem{thm}{\bf Theorem}
\begin{thm} (Relation between stochastic governing law and Fokker-Planck equation)\\
\label{thm1}
For every $\varepsilon >0$, the probability density function $p\left( x,t|z,0 \right)$ and the jump measure, drift and diffusion have the following relations:\\
1) For every $x$ and $z$ satisfying $\left| x-z \right|>\varepsilon $,
\begin{align}\label{T.1}
\underset{t\to 0}{\mathop{\lim }}\,{{t}^{-1}}{p}\left( x,t|z,0 \right)=\sigma^{-n}W\left( \sigma^{-1}\left( x-z \right) \right)
\end{align}
uniformly in $x$ and $z$.\\
2) For $i=1,\ 2,\ \ldots ,\ n$,
\begin{align}\label{T.2}
\underset{t\to 0}{\mathop{\lim }}\,{{t}^{-1}}\int_{\left| x-z \right|<\varepsilon}{\left( {{x}_{i}}-{{z}_{i}} \right)p\left( x,t|z,0 \right) \textrm{d}\mathbf{x}}={{b}_{i}}\left( z \right).
\end{align}
3) For $i,j=1,\ 2,\ \ldots ,\ n$,
\begin{align}\label{T.3}
\underset{t\to 0}{\mathop{\lim }}\,{{t}^{-1}}\int_{\left| x-z \right|<\varepsilon}{\left( {{x}_{i}}-{{z}_{i}} \right)\left( {{x}_{j}}-{{z}_{j}} \right)p\left( x,t|z,0 \right) \textrm{d}\mathbf{x}}={{a}_{ij}}\left( z \right)+\sigma^{-n}\int_{\left|y\right|<\varepsilon}{{y}_{i}{y}_{j}W\left( {\sigma}^{-1}{y} \right) dy}.
\end{align}
\end{thm}

This theorem can be reformulated as the following corollary, the so-called nonlocal Kramers-Moyal formulas, which express the jump measure, drift and diffusion via the sample paths of the stochastic differential equation (\ref{SDE}).

\newtheorem{cor}[thm]{\bf Corollary}
\begin{cor}(Nonlocal Kramers-Moyal formulas)\\
\label{cor2}
For every $\varepsilon >0$, the sample path solution $x\left( t \right)$ of the stochastic differential equation (\ref{SDE}) and the jump measure, drift and diffusion have the following relations:\\
1) For every $m>1$,
\begin{align}\label{C.1}
\underset{t\to 0}{\mathop{\lim }}\,{{t}^{-1}}\mathbb{P}\left\{ \left. \left| {x}\left( t \right)-{z}\right| \in \left[ \varepsilon,\ m\varepsilon \right) \right| x\left( 0 \right)=z \right\}=\sigma ^{-n} \int_{ \left| y \right| \in \left[ \varepsilon,\ m\varepsilon \right)} {W\left( \sigma ^{-1} {y} \right) dy}.
\end{align}
2) For $i=1,\ 2,\ \ldots ,\ n$,
\begin{align}\label{C.2}
\underset{t\to 0}{\mathop{\lim }}\,{{t}^{-1}}\mathbb{P}\left\{ \left. \left| x\left( t \right)-z \right| <\varepsilon  \right| x\left( 0 \right) =z \right\}\cdot \mathbb{E} \left[ \left. \left( {{x}_{i}}\left( t \right) -{{z}_{i}} \right) \right| x\left( 0 \right)=z;\ \left| x\left( t \right)-z \right| <\varepsilon  \right] ={{b}_{i}}\left( z \right).
\end{align}
3) For $i,j=1,\ 2,\ \ldots ,\ n$,
\begin{align}\label{C.3}
  & \underset{t\to 0}{\mathop{\lim }}\,{{t}^{-1}}\mathbb{P}\left\{ \left. \left| x\left( t \right)-z \right| <\varepsilon  \right| x\left( 0 \right) =z \right\}\cdot \mathbb{E} \left[ \left. \left( {{x}_{i}}\left( t \right) -{{z}_{i}} \right) \left( {{x}_{j}}\left( t \right) -{{z}_{j}} \right) \right| x\left( 0 \right)=z;\ \left| x\left( t \right)-z \right| <\varepsilon  \right] \nonumber\\
 & ={{a}_{ij}}\left( z \right)+\sigma^{-n}\int_{\left|y\right|<\varepsilon}{{y}_{i}{y}_{j}W\left( {\sigma}^{-1}{y} \right) dy}.
\end{align}
\end{cor}

In this work, we consider a special but significant L\'evy motion for the sake of concreteness due to its extensive physical applications  \cite{OEBbook} as an example to illustrate our method,  the rotational symmetric  $\alpha$-stable L\'evy motion \cite{LMbook}. Its detailed information is present in the Appendix. In this case, the identification of the L\'evy jump measure is transformed to learn the stability parameter $\alpha$ from data.

We have shown that the relationship between the nonlocal Kramers-Moyal formulas and the coefficients of stochastic differential equation. The nonlocal Kramers-Moyal formulas provide a way to estimate the L\'evy jump measure, drift coefficient and diffusion coefficient using transition probability density. This is naturally associated with normalizing flows.

\section{Method}
In this section, we will state how to learn the coefficients of stochastic differential equation by combining the nonlocal Kramers-Moyal formulas with normalizing flows. To clarify which formulas in section \ref{KMformulas} we used to estimate the coefficients, we provide a practical implementation here.\\

\textbf{Implementation} The practical implementation for learning the coefficients of stochastic differential equation:
\begin{itemize}
\item [1] Estimating transition probability density $p(x,t|z,0)$ from short bursts of simulation data with inital value $z$ using neural spline flows (\ref{NSF_trans}) or RealNVP (\ref{RealNVP_Trans}).
\item [2] Using the inverse transformation of normalizing flows $T_{\theta}^{-1}$ to obtain the samples from the density estimated $p(x,t|z,0)$, and then the formula (\ref{C.1}) is used to estimate the L\'evy jump measure, i.e., anomalous diffusion coefficient $\sigma$ and stability parameter $\alpha$. Here we use the resampled samples to estimate the probability on the left-hand side of the formula (\ref{C.1}). For more details, see \cite{YangLi2020a}.
\item[3] Using the transformation of normalizing flows $T_{\theta}$ to estimate the density $p(x,t|z,0)$, then we can obtain the drift coefficient and diffusion coefficient from (\ref{T.2}) and (\ref{T.3}).
\end{itemize}

\section{Results}\label{results}
In the previous sections, we have introduced the nonlocal Kramers-Moyal formulas and normalizing flows. The nonlocal Kramers-Moyal formulas show the relationship between the transition probability density $p(x,t|z,0)$ and coefficients (L\'evy measure, drift and diffusion) for corresponding stochastic differential equation. Moreover, the normalizing flows provide a flexible way to approximate the transition density $p(x,t|z,0)$ from sample paths data. Therefore, we can estimate the coefficients of stochastic differential equation from sample paths.\\

Let $t^*$ be a value of time "near" zero, a dataset $\mathcal{D}=\{x_{i}^{(z)}\}_{i=1}^{n}$ sampled from $p_{x}(x,t^*|z,0)$. To be specific, these samples come from simulating stochastic differential equation with initial value $z$ using Euler-Maruyama method, where $z\in[-2,2]$ or $z\in[-2,2]\times[-2,2]$.\\

\subsection{A one-dimensional system}
\begin{figure}
\centering

\includegraphics[trim={0cm 0cm 0cm 0cm},clip,width=\textwidth]{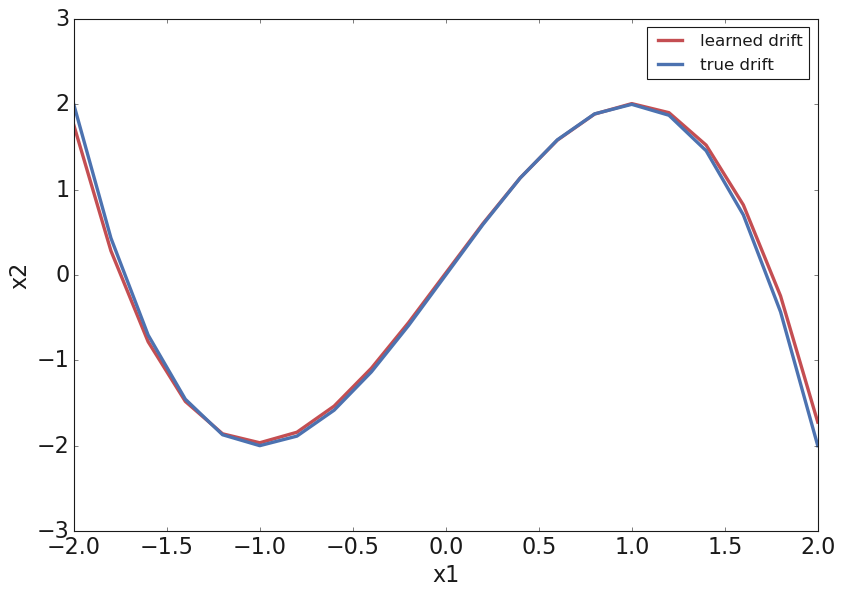}
\caption{1D system: Learning the Drift coefficient. The blue line is the true drift coefficient, the red line is the learned drift coefficient.}
\label{1D_drift}
\end{figure}

Consider a one-dimensional stochastic differential equation with pure jump L\'evy motion
\begin{align}\label{EX1}
&dX(t)=(3X(t)-X^3(t))dt+dL^{\alpha}(t),
\end{align}
where $L_{1}^{\alpha}$ is a scalar, real-value $\alpha-$stable L\'evy motion with triple $(0,0,\nu_{\alpha})$. Here we take $\alpha=1.5$. The jump measure is given by
\begin{align}
\nu_{\alpha}(dx)=c(1,\alpha)||y||^{-1-\alpha}dx,
\end{align}
with $c(1,\alpha)=\frac{\alpha\Gamma((1+\alpha)/2)}{2^{1-\alpha}\pi^{\frac{1}{2}}\Gamma(1-\alpha/2)}$. Here we take $t^*=0.01$, sample size $n=5000$ and standard normal density as our prior density. For the transformation of normalizing flows, we take neural spline flows as our basis transformation, denoted by $T_{\theta}$. In order to improve the complexity of transformation, we take the $N=32$ compositions $T_{\theta}^N$ as our final transformation, where the neural network in (\ref{NSF_trans}) is a full connected neural network with 3 hidden layers and 32 nodes each layer. In addition, we choose the interval $[-B, B]=[-3,3]$ and $K=5$ bins in section \ref{NSF}.\\
We compare the true drift coefficient and the learned drift coefficient in Figure \ref{1D_drift}. In addition, the estimated $\alpha=1.55, \sigma=1.09$, where the true values $\alpha=1.5, \sigma=1$.\\

In this subsection, we showed how to estimate the coefficients of stochastic differential equation using neural spline flows and nonlocal Kramers-Moyal formulas. However, the biggest disadvantage of neural spline flows is that the root of rational polynomial is required when solving the inverse transformation $T_{\theta}$. Therefore we prefer to use another method to approximate multidimensional probability density, i.e., RealNVP here.

As we know in Section \ref{RealNVP}, the transformation (\ref{RealNVP_Trans}) is our basis transformation. For the target probability density $p_{x}(x),x\in \mathbb{R}^2$, our final transformation consists of compositions of the following three types of prototypical transformations:
\begin{align}
    z_{1}&=x_{1}, \nonumber \\
    z_{2}&=\frac{1}{C}[x_{2} e^{\mu_{1}(x_{1})}+\nu_{1}(x_{1})],
\end{align}

\begin{align}
    z_{1}&=\frac{1}{C}[x_{1} e^{\mu_{2}(x_{2})}+\nu_{2}(x_{2})], \nonumber \\
    z_{2}&=x_{2},
\end{align}

\begin{align}
    z_{1}&=x_{1}, \nonumber \\
    z_{2}&=\frac{1}{C}[x_{2} e^{\mu_{3}(x_{1})}+\nu_{3}(x_{1})],
\end{align}
where $(\mu_i,\nu_i),i=1,2,3$ are six different full connected neural networks with 3 hidden layers and 16 nodes each layer. Denoted this final transformation by $T_{\theta}$.\\
\subsection{A decoupled system}\label{example2}

Consider a two-dimensional stochastic differential equation with a pure jump L\'evy motion
\begin{align}\label{EX2}
&dX_{1}(t)=(3X_{1}(t)-X_{1}^3(t))dt+dL_{1}^{\alpha}(t), \nonumber \\
&dX_{2}(t)=(3X_{2}(t)-X_{2}^3(t))dt+dL_{2}^{\alpha}(t),\\
 \nonumber
\end{align}
where $(L_{1}^{\alpha},L_{2}^{\alpha})$ is a two-dimensional rotational symmetry real-value $\alpha-$stable L\'evy motions with triple $(0,0,\nu_{\alpha})$. Here we take $\alpha=1.5$. The jump measure is given by
\begin{align}
\nu_{\alpha}(dx)=c(2,\alpha)||y||^{-2-\alpha}dx,
\end{align}
with $c(2,\alpha)=\frac{\alpha\Gamma((2+\alpha)/2)}{2^{1-\alpha}\pi^{\frac{1}{2}}\Gamma(1-\alpha/2)}$. Here we take $t^*=0.05$, sample size $n=10000$, hyperparameter $C=\frac{1}{3}$ and standard normal density as our prior density.

\begin{figure}
\centering
\includegraphics[trim={0cm 0cm 0cm 0cm},clip,width=\textwidth]{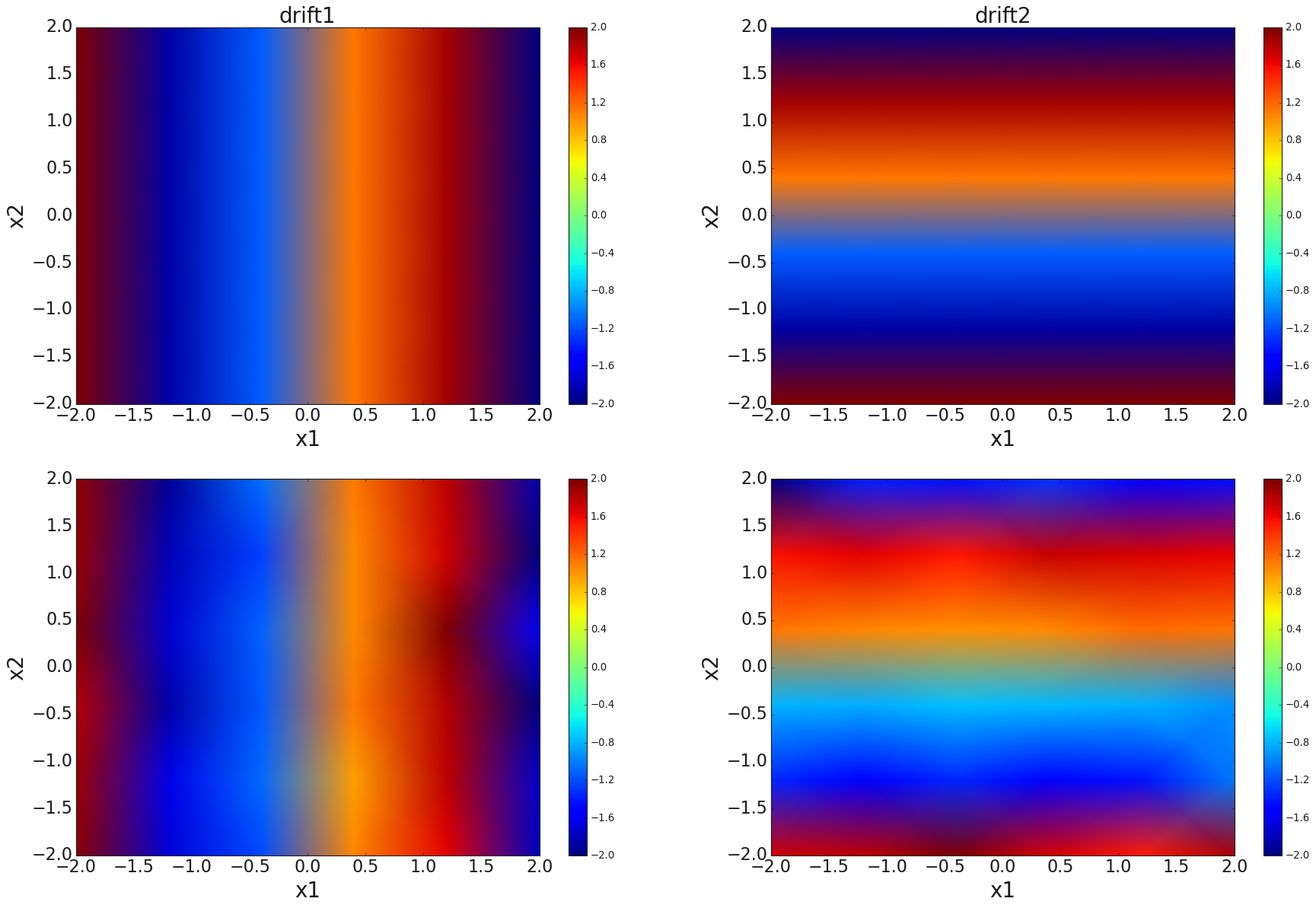}
\caption{2D decoupled system with pure jump L\'evy motion: Learning the drift coefficients. Top row: The true values of drift coefficients. Bottom row: The learned values of drift coefficients.}
\label{2D_drift_decoupling}
\end{figure}

We compare the true drift coefficients and the learned drift coefficients in Figure \ref{2D_drift_decoupling}.  In addition, the estimated $\alpha=1.41, \sigma=1.04$, where the true values $\alpha=1.5, \sigma=1$.\\

Example \ref{example2} showed us we can learn the drift term of such a decoupling stochastic system from data. Next, we will illustrate our method for a more complexity example.\\

\subsection{A coupled system}
Consider a two-dimensional coupling stochastic differential equation with multiplicative Brownian motion and  L\'evy motion
\begin{align}\label{EX3}
&dX_{1}(t)=(0.001X_{1}(t)-X_{1}(t)X_{2}(t))dt+X_{1}(t)dB_{1}(t)+dL_{1}^{\alpha}(t), \nonumber \\
&dX_{2}(t)=(-6X_{2}(t)+0.25X_{1}^2(t))dt+X_{2}(t)dB_{2}(t)+dL_{2}^{\alpha}(t),\\
\end{align}
where $(B_{1},B_{2})$ is a two-dimensional Brownian motion, $(L_{1}^{\alpha},L_{2}^{\alpha})$ is a two-dimensional rotational symmetry real-value $\alpha-$stable L\'evy motions with triple $(0,0,\nu_{\alpha})$. Here we take $\alpha=1.5$.  The jump measure is given by
\begin{align}
\nu_{\alpha}(dx)=c(2,\alpha)||y||^{-2-\alpha}dx,
\end{align}
with $c(2,\alpha)=\frac{\alpha\Gamma((2+\alpha)/2)}{2^{1-\alpha}\pi^{\frac{1}{2}}\Gamma(1-\alpha/2)}$. Here we take $t^*=0.01$, sample size $n=10000$, hyperparameter $C=\frac{1}{3}$ and standard normal density as our prior density.

\begin{figure}
\centering

\includegraphics[trim={0cm 0cm 0cm 0cm},clip,width=\textwidth]{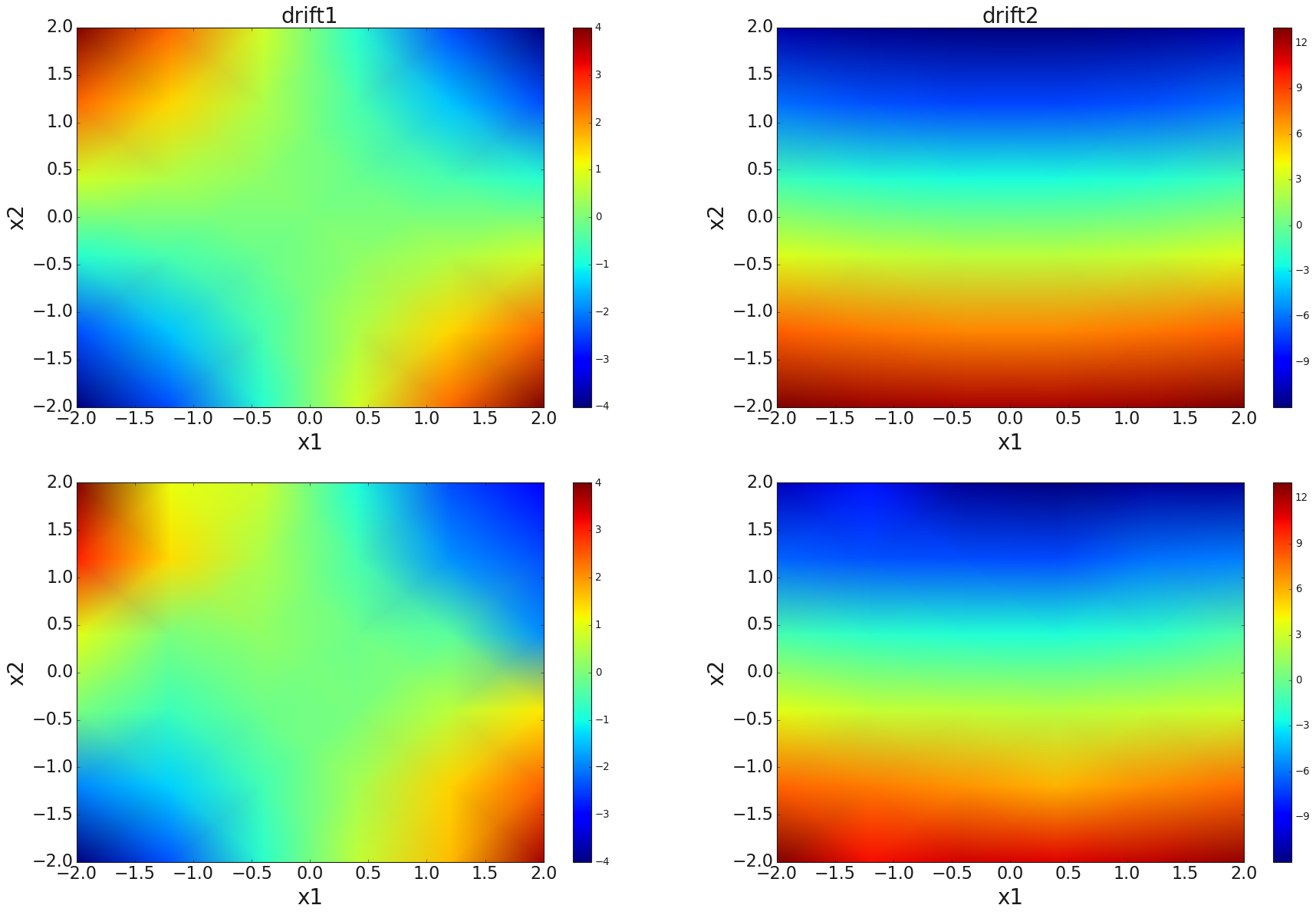}
\caption{2D coupled system with multiplicative Brownian motion and L\'evy motion: Learning the drift coefficients. Top row: The true values of drift coefficients. Bottom row: The learned values of drift coefficients.}
\label{2D_drift_coupling}
\end{figure}

\begin{figure}
\centering
\includegraphics[trim={0cm 0cm 0cm 0cm},clip,width=\textwidth]{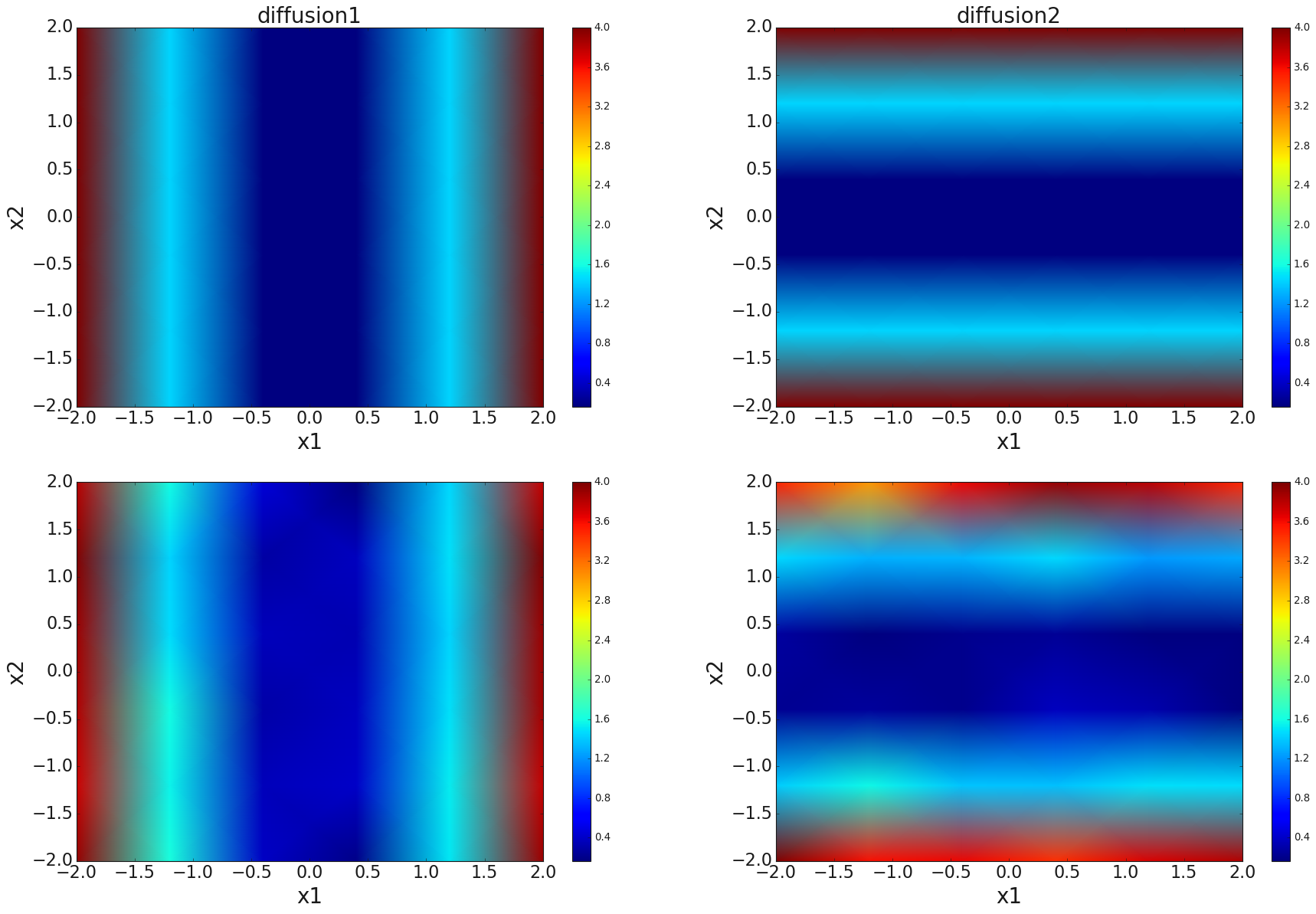}
\caption{2D coupled system with multiplicative Brownian motion and L\'evy motion: Learning the diffusion coefficients of Brownian motion. Top row: The true values of diffusion coefficients. Bottom row: The learned values of diffusion coefficients. }
\label{2D_diffusion_coupling}
\end{figure}

We compare the true drift coefficients and the learned drift coefficients in Figure \ref{2D_drift_coupling}. Figure \ref{2D_diffusion_coupling} show the learned diffusion coefficients.  In addition, the estimated $\alpha=1.58, \sigma=1.23$, where the true values $\alpha=1.5, \sigma=1$. \\

\section{Discussion}
In this work, we have presented a new method to extract stochastic governing laws by nonlocal Kramers-Moyal formulas. In particular, we use normalizing flows to estimate transition probability density from sample path data, and then identify the L\'evy jump measure, drift coefficient and diffusion coefficient of a stochastic differential equation. We further develop a numerical algorithm and test it on some prototypical examples to confirm its efficacy.

Comparing with the original data-driven algorithm \cite{YangLi2020a}, our method reduces the amount of data and does not require the basis to represent these coefficients in the stochastic governing laws. This makes it easier to extract stochastic governing laws from noisy data with non-Gaussian statistical features.  Although we have improved our previous work, there are still some challenges. We reduced the amount of data, but the accuracy of our estimation results was slightly worse than the original one. That is because we did not express the coefficients in terms of a basis (whose selection is also a challenge). The accuracy of the estimate also depends on individual dynamical systems. How to design an algorithm with strong robustness and high estimation accuracy is our future work.
\section*{Acknowledgments}
The authors are grateful to Ting Gao and Cheng Fang for helpful discussions.

\section*{Data Availability}
 Research code is shared via \url{https://github.com/Yubin-Lu/Extracting-Stochastic-Governing-Laws-by-Nonlocal-Kramers-Moyal-Formulas}.


\setcounter{equation}{0}
\renewcommand\theequation{A.\arabic{equation}}
\section*{Appendix}\label{App}
\subsection*{L\'evy motions}
Let $L=(L_{t},t\geq0)$ be a stochastic process defined on a probability space $(\Omega, \mathcal{F}, P)$. We say that L is a L\'evy motion if:
(1) $L_{0}=0 (a.s.)$;\\
(2) $L$ has independent and stationary increments;\\
(3) $L$ is stochastically continuous, i.e., for all $a>0$ and for all $s\geq0$
$$
\lim\limits_{t\to s}P(|L_{t}-L_{s}|>a)=0.
$$

\subsection*{Characteristic functions}
For a L\'evy motion $(L_{t},t\geq0)$, we have the L\'evy-Khinchine formula,
$$
\mathbb{E}[e^{i(u,L_t)}] = exp\{t[i(b,u)-\frac{1}{2}(u,Au)+\int_{\mathbb{R}^d\backslash\{0\}} [e^{i(u,y)}-1-i(u,y)I_{\{\Vert y\Vert<1\}}(y)]\, \nu(dy)]\},
$$
for each $t\geq0$, $u\in\mathbb{R}^n$, where $(b,A,\nu)$ is the triple of L\'evy motion $(L_{t}, t\geq0)$.\\

\textbf{Theorem (The L\'evy-It\^{o} decomposition)} If $(L_{t},t\geq0)$ is a L\'evy motion with $(b,A,\nu)$, then there exists $b\in\mathbb{R}^n$, a Brownian motion $B_{A}$ with covariance matrix $A$ and an independent Poisson random measure $N$ on $\mathbb{R}^{+}\times(\mathbb{R}^n-\{0\})$ such that, for each $t\geq0$, \\
$$
L_{t}=bt+B_{A}(t)+\int_{0<|x|<c} x\, \tilde{N}(t,dx)+\int_{|x|\geq c} x\, N(t,dx),
$$
where $\int_{|x|\geq c} x\, N(t,dx)$ is a Poisson integral and $\int_{0<|x|<c} x\, \tilde{N}(t,dx)$ is a compensated Poisson integral defined by
$$
\int_{0<|x|<c} x\, \tilde{N}(t,dx)=\int_{0<|x|<c} x\, N(t,dx)-t\int_{0<|x|<c} x\, \nu(dx).
$$
\\

\subsection*{The $\alpha$-stable L\'evy motions}
The $\alpha$-stable L\'evy motion is a special but most popular type of the L\'evy process defined by the stable random variable with the distribution ${{S}_{\alpha }}\left( \delta ,\ \beta ,\ \lambda  \right)$ \cite{HeavytailBook, LMBook2,OEBbook}. Usually, $\alpha \in \left( 0,\ 2 \right]$ is called the stability parameter, $\delta \in \left( 0,\ \infty  \right)$ is the scaling parameter, $\beta \in \left[ -1,\ 1 \right]$ is the skewness parameter and $\lambda \in \left( -\infty ,\ \infty  \right)$ is the shift parameter.

A stable random variable $X$ with $0<\alpha <2$ has the following "heavy tail" estimate:
$$
\underset{x\to \infty }{\mathop{\lim }}\,{{y}^{\alpha }}\mathbb{P}\left( X>y \right)={{C}_{\alpha }}\frac{1+\beta }{2}{{\delta }^{\alpha }},
$$
where ${{C}_{\alpha }}$ is a positive constant depending on $\alpha$. In other words, the tail estimate decays polynomially. The $\alpha$-stable L\'evy motion has larger jumps with lower jump frequencies for smaller $\alpha$ ($0<\alpha <1$), while it has smaller jump sizes with higher jump frequencies for larger $\alpha$ ($1<\alpha <2$). The special case $\alpha=2$ corresponds to (Gaussian) Brownian motion. For more information aboutL\'evy process, refer to References \cite{Duan2015, Applebaum}.

\subsection*{Fokker-Planck equations}
Consider a stochastic differential equation in $\mathbb{R}^n$:
\begin{align}\label{sde1}
dX_t=f(X_{t-})dt+\sigma(X_{t-})dB_t + dL_{t}^{\alpha},\quad X_0=x_0,
\end{align}
where $f$ is a vector field, $\sigma$ is a $n\times n$ matrix, $B_t$ is Brownian motion in $\mathbb{R^n}$ and $L_{t}^{\alpha}$ is a symmetric $\alpha-$stable L\'evy moton in $\mathbb{R^n}$ with the generating triplet $(0, 0, \nu_{\alpha})$. The jump measure
\begin{align}
\nu_{\alpha}=c(n,\alpha)\parallel y\parallel^{-n+\alpha}dy,
\end{align}
with $c(n,\alpha)=\frac{\alpha\Gamma((n+\alpha)/ 2)}{2^{1-\alpha}\pi^{n/2}\Gamma(1-\alpha/ 2)}$. The processes $B_t$ and $L_t^\alpha$ are taken to be independent.\\
The Fokker-Planck equation for the stochastic differential equation  (\ref{sde1}) is then
\begin{align}\label{FPeqn}
p_t&=-\bigtriangledown\cdot(fp)+\frac{1}{2}{\rm Tr}[H(\sigma\sigma^{T}p)]\nonumber\\
&+\int_{\mathbb{R}^n\backslash\{0\}} [p(x+y,t)-p(x,t)]\, \nu_{\alpha}(dy),
\end{align}
where $p(x,0)=\delta(x-x_0)$. See \cite{Applebaum,Duan2015} for more details.\\

\end{document}